\documentclass[10pt,a4paper]{article}
\usepackage[a4paper]{geometry}
\usepackage{amssymb,latexsym,amsmath,amsfonts,amsthm}
\usepackage{graphicx}

\usepackage{epsfig}

\newtheorem{theorem}{Theorem}[section]
\newtheorem{lemma}[theorem]{Lemma}

\theoremstyle{definition}
\newtheorem{definition}[theorem]{Definition}

\theoremstyle{remark}

\newtheorem{remark}[theorem]{Remark}

\numberwithin{equation}{section}

\title{Hermite-Pad\'e approximation for certain systems of meromorphic functions}

\date{\today}

\author{U Fidalgo Prieto\footnotemark[1], G L\'opez Lagomasino\footnotemark[2], S Medina Peralta\footnotemark[2]}

\begin{document}

\maketitle

\renewcommand{\thefootnote}{\fnsymbol{footnote}}
\footnotetext[1]{Department of Mathematics, University of Mississippi, Hume Hall 327, Oxford,  MS 38677-1848, USA.}
\footnotetext[2]{Departamento de Matem\'{a}ticas, Universidad Carlos III de Madrid, Avda. Universidad
30, 28911 Legan\'{e}s, Madrid, Spain. email: lago\symbol{'100}math.uc3m.es.}
\footnotetext[3]{The authors received  support from research grant MTM2012-36372-C03-01 of Ministerio de Econom\'{\i}a y Competitividad, Spain.}

\begin{abstract}
We study the convergence of sequences of type I and type II Hermite-Pad\'e approximants for certain  systems of meromorphic functions made up of  rational modifications of Nikishin systems of functions.
\end{abstract}



\maketitle

\section{Introduction}

This article deals with the convergence of sequences of type I and type II Hermite-Pad\'e approximants. In the study of the algebraic and analytic properties of such approximants the contributions of A.A. Gonchar and H. Stahl have been pioneering and significant, see a.e.  \cite{DrSt0}, \cite{DrSt1}, \cite{DrSt2}, \cite{GR},   \cite{GRS}, and \cite{Stahl}.

Since their introduction in \cite{Her} by Ch. Hermite, these approximants have been an invaluable tool in number theory (for a survey of such results see \cite{walter}). They may be used in the construction of simultaneous quadrature rules, see for example \cite{cou}  and \cite{LIF2}. Recently, they have appeared in models coming from random matrix theory and mathematical physics  (a summary of such models is contained in \cite{Kuij}).

\begin{definition}
Let ${\bf f}=(f_1,\ldots,f_m)$ be a system of $m$ formal power series where
\begin{equation}\label{series}
f_j(z)= c_{0,j} + \frac{c_{1,j}}{z}+\frac{c_{2,j}}{z^2}+\cdots, \qquad j=1,\ldots,m.
\end{equation}
Fix a multi-index ${\bf n}=(n_1,\ldots,n_m)\in {\mathbb{Z}}_+^m \setminus \{{\bf 0}\}$, ${\mathbb{Z}}_+=\{0,1,2,\ldots\}.$  The vector rational function  $\displaystyle {\bf R}_{\bf n}=\left(\frac{P_{{\bf n},1}}{Q_{\bf n}},\ldots,\frac{P_{{\bf n},m}}{Q_{\bf n}}\right)$ is called a type II Hermite-Pad\'e approximant with respect to ${\bf f}$ and ${\bf n}$ if the polynomials $Q_{\bf n}$, $P_{{\bf n},j},$ $j=1,\ldots,m$ satisfy:
\begin{itemize}
\item[i)]  $\deg Q_{\bf n} \leq n_1+\cdots+n_m$, $Q_{\bf n}\not \equiv 0$,
\item[ii)] $
(Q_{\bf n}  f_j -P_{{\bf n},j})(z)=\mathcal{O}(1/z^{n_j+1}),  \qquad j=1,\ldots,m.$
\end{itemize}
Similarly, the vector polynomial $(a_{{\bf n},0},\ldots,a_{{\bf n},m})$ is  a type I Hermite-Pad\'e approximant of $\bf f$ with respect to $\bf n$ if:
\begin{itemize}
\item[iii)] $\deg a_{{\bf n},j} \leq n_j -1, j=1,\ldots,m,$ not all identically equal to zero,
\item[iv)] $(a_{{\bf n},0} + \sum_{j=1}^m a_{{\bf n},j}f_j)(z) = \mathcal{O}(1/z^{n_1+\cdots+n_m}).$
\end{itemize}
The right hand of ii) and iv) represent formal expansions in descending powers of $z$ and equality is in the sense of the coefficients.
\end{definition}

Given $\bf f$ and $\bf n$, it is easy to see that  type I and type II Hermite-Pad\'e approximation always exist. Their construction reduces to solving  homogeneous linear systems  with one equation less  than unknowns.  When $m=1$ both definitions coincide with  classical diagonal Pad\'e approximation. In contrast with Pad\'e approximation, when $m \geq 2$ the uniqueness of these approximants is not guaranteed  in general.

We will focus our study on an important class of systems of functions. Due to classical results of A.A. Markov \cite{Mar} and T. Stieltjes \cite{St}, the Cauchy transform of measures supported on the real line $\mathbb{R}$ play a central role  in the convergence theory of diagonal Pad\'e approximation. For type I and type II Hermite Pad\'e approximants, to a large extent a similar place is occupied by the so called Nikishin systems of functions introduced in \cite{Nik}. In particular, for such systems Markov and Stieltjes type theorems are verified (see, for example, \cite{Bus},  \cite{LF3}, \cite{GRS}, \cite{LS} and \cite{Stahl}).

In the sequel $\Delta$ denotes an interval contained in a half line of the real axis. By $\mathcal{M}(\Delta)$ we denote the class of all Borel measures $s$ with constant sign whose support, contained in $\Delta$, has infinitely many points, and $x^\nu \in L_1(s)$ for all $\nu \in \mathbb{Z}_+$. Set
\[ \displaystyle \widehat{s}(z) = \int\frac{d s(x)}{z-x}.
\]
Then
\begin{equation}\label{expansion}
\displaystyle \widehat{s}(z)  \sim  \sum_{j=0}^{\infty} \frac{c_j}{z^{j+1}}, \qquad c_j = \int x^j ds(x).
\end{equation}
If  the support  of $s$, $\mbox{supp}(s)$, is bounded  the series is convergent in a neighborhood of $\infty$; otherwise, the expansion is asymptotic at $\infty$. That is, for each $k \geq 0$
\[ \lim_{z \to \infty} z^{k+1}\left(\widehat{s}(z) - \sum_{j=0}^{k-1} \frac{c_j}{z^{j+1}}\right) = c_k,
\]
where the limit is taken along any curve which is non tangential to $\mbox{supp}(s)$ at $\infty$.

Let $\Delta \subset \mathbb{R}_+$. Stieltjes' theorem (which contains Markov's theorem when $\Delta$ is bounded) states that if the moment problem for $(c_n)_{n \in \mathbb{Z}_+}$ is determinate then
\begin{equation} \label{MS} \lim_{n\to \infty} \frac{P_n}{Q_n}(z) = \widehat{s}(z),
\end{equation}
uniformly on each compact subset of (inside) $\mathbb{C} \setminus \Delta$. Here, $\{P_n/Q_n\}_{n\geq 0}$ denotes the diagonal sequence of Pad\'e approximants of $\widehat{s}$. A sufficient condition, due to T. Carleman \cite{Car}, for the moment problem to be determinate is
\begin{equation}\label{Carle} \sum_{n \geq 0} |c_n|^{-1/2n} = \infty,
\end{equation}

In an attempt to extend Markov's theorem to a general class of meromorphic functions, A.A. Gonchar considered functions of the form $\widehat{s} + r$ where $r$ is a rational function whose poles lie in $\mathbb{C} \setminus \Delta$. In \cite{Gon0}, he proved that if  $\Delta$ is a bounded interval and the sequence of orthogonal polynomials with respect to $s$ have ratio asymptotic, then  (\ref{MS}) takes place showing, additionally, each pole of $r$ in $\mathbb{C} \setminus \Delta$ ``attracts'' as many zeros of $Q_n$ as its order and the remaining zeros of $Q_n$ accumulate on $\Delta$ as $n \to \infty$. Later, in \cite{Rak1} E.A. Rakhmanov obtained a full extension of Markov's theorem when $r$ has real coefficients and proved that if $r$ has complex coefficients then such a result is not possible  without extra assumptions on the measure $s$. The case of unbounded $\Delta$ was solved in \cite{L1}, when $r$ has real coefficients, and \cite{L2}, when $r$ has complex coefficients.

Let us introduce Nikishin systems. Let $\Delta_{\alpha}, \Delta_{\beta}$ be two intervals contained in the real line with at most a common point. Take $\sigma_{\alpha} \in {\mathcal{M}}(\Delta_{\alpha})$ and $\sigma_{\beta} \in {\mathcal{M}}(\Delta_{\beta})$ such that $\widehat{\sigma}_{\beta} \in L_1(\sigma_{\alpha})$. Using the  differential notation, we define a third measure $\langle \sigma_{\alpha},\sigma_{\beta} \rangle$ as follows
\[d \langle \sigma_{\alpha},\sigma_{\beta} \rangle (x) := \widehat{\sigma}_{\beta}(x) d\sigma_{\alpha}(x).\]
In consecutive products of measures such as $\langle \sigma_{\gamma},  \sigma_{\alpha},\sigma_{\beta} \rangle :=\langle \sigma_{\gamma}, \langle \sigma_{\alpha},\sigma_{\beta} \rangle \rangle,$ we assume not only that $\widehat{\sigma}_{\beta} \in L_1(\sigma_{\alpha})$ but also $\langle \sigma_{\alpha},\sigma_{\beta} \widehat{\rangle} \in L_1(\sigma_{\gamma})$, where $\langle \sigma_{\alpha},\sigma_{\beta} \widehat{\rangle}$ denotes the Cauchy transform of $\langle \sigma_{\alpha},\sigma_{\beta}  {\rangle}$.

\begin{definition}
Take a collection  $\Delta_j, j=1,\ldots,m,$ of intervals such that
\[ \Delta_j \cap \Delta_{j+1} = \emptyset, \qquad \mbox{or} \quad \Delta_j \cap \Delta_{j+1} = \{x_{j,j+1}\}, \quad j=1,\ldots,m-1,
\]
where $x_{j,j+1}$ is a single point. Let $(\sigma_1,\ldots,\sigma_m)$ be a system of measures such that $\mbox{\rm Co}(\mbox{\rm supp} (\sigma_j)) = \Delta_j, \sigma_j \in {\mathcal{M}}(\Delta_j), j=1,\ldots,m,$  and
\begin{equation}\label{eq:autom}
\langle \sigma_{j},\ldots,\sigma_k  {\rangle} := \langle \sigma_j,\langle \sigma_{j+1},\ldots,\sigma_k\rangle\rangle\in {\mathcal{M}}(\Delta_j),  \qquad  1 \leq j < k\leq m,
\end{equation}
where $\mbox{\rm Co}(E)$ denotes the convex hull of the set $E$.
When $\Delta_j \cap \Delta_{j+1} = \{x_{j,j+1}\}$ we also assume that $x_{j,j+1}$ is not a mass point of either $\sigma_j$ or $\sigma_{j+1}$.
We say that ${\bf s}=(s_{1,1},\ldots,s_{1,m}) = {\mathcal{N}}(\sigma_1,\ldots,\sigma_m)$, where
\[ s_{1,1} = \sigma_1, \quad s_{1,2} = \langle \sigma_1,\sigma_2 \rangle, \ldots \quad , \quad s_{1,m} = \langle \sigma_1, \sigma_2,\ldots,\sigma_m  \rangle,
\]
is the Nikishin system of measures generated by $(\sigma_1,\ldots,\sigma_m)$. The corresponding Nikishin system of functions will be denoted  by ${\bf \widehat{s}}=\left(\widehat{s}_{1,1},\ldots,\widehat{s}_{1,m}\right)$, where $\widehat{s}_{1,j}$ is the Cauchy transform of $s_{1,j}$.
\end{definition}

This definition extends the one given in  \cite{Nik} by allowing the generating measures to have unbounded support and/or have consecutive $\Delta_j$ with a common endpoint.   In what follows for $1\leq j\leq k\leq m$ we denote
\[ s_{j,k} := \langle \sigma_j,\sigma_{j+1},\ldots,\sigma_k \rangle, \qquad s_{k,j} := \langle \sigma_k,\sigma_{k-1},\ldots,\sigma_j \rangle.
\]

Observe that the functions $\widehat{s}_{1,j}$ are holomorphic in the complement of $\Delta_1$; that is,  $\widehat{s}_{1,j}\in {\mathcal{H}}\left( {\mathbb{C}}\setminus \Delta_1\right)$, $j=1,\ldots,m$. Consider a vector of rational functions ${\bf r}=(r_1,\ldots,r_m)=\left(\displaystyle \frac{v_{1}}{t_{1}},\ldots, \frac{v_{m}}{t_{m}}\right)$, such that $\deg t_{j}=d_{j}$ and  $\deg v_{j}<d_{j}$, for every $j=1,\ldots,m$. We assume that $v_j/t_j, j=1,\ldots,m$ is irreducible.  We consider systems of meromorphic functions  of the form  ${\bf f}=(f_{1},\ldots, f_{m}) = \widehat{\bf s} + {\bf r}$, where
\begin{equation}\label{fs}
f_{j}(z)=\widehat{s}_{1,j}(z)+ r_{j}(z), \qquad j=1,\ldots,m.
\end{equation}

A  Stieltjes type theorem was proved in \cite{Bus} for type II Hermite-Pad\'e aproximants of a Nikishin system $\widehat{\bf s}$. Our goal is to provide an analogue for systems of functions of the form ${\bf f} = \widehat{\bf s} + {\bf r}$. For $m=2$, a partial result appears in \cite{Bus2}. In the sequel, we assume that $m\geq 2$ since the case $m=1$ corresponds to Pad\'e approximation and there the corresponding results are known.

Let $\mbox{\rm dist}(E,F)$ denote the distance in the spherical metric of $\mathbb{C}$ between the sets $E$ and $F$. Let $\mathcal{C}(\Delta)$ be the class of all measures in $\mathcal{M}(\Delta)$ such that after an affine transformation which takes $\Delta$ to a subset of $\mathbb{R}_+$ the image measure satisfies (\ref{Carle}) (of course the result does not depend on the affine transformation taken). We have

\begin{theorem}\label{TCTII}
Let ${\bf s} = (s_{1,1},\ldots,s_{1,m}) = \mathcal{N}(\sigma_1,\ldots,\sigma_m)$ and ${\bf r} = (r_1,\ldots,r_m)$ be given,  where the rational functions $r_j, j=1,\ldots,m$ have real coefficients, for different $j$ their poles are distinct, and they all lie in  ${\mathbb{C}}\setminus (\Delta_1 \cup \Delta_m)$. Let $\Lambda \subset \mathbb{Z}_+^{m}$ be an infinite sequence of distinct multi-indices  such that
\begin{equation}\label{cond1} \sup \left\{\max_{j=1,\ldots,m}(n_j) - \min_{k=1,\ldots,m}(n_k) : {\bf n} =(n_1,\ldots,n_m) \in \Lambda\right\} \leq C < \infty. \end{equation}
Assume that $\mbox{\rm dist}(\Delta_{1},\Delta_2) >0$  or $\sigma_1 \in \mathcal{C}(\Delta_1)$  , and $\mbox{\rm dist}(\Delta_{m-1},\Delta_m) >0$  or $\sigma_m \in \mathcal{C}(\Delta_m)$. Let $({\bf R}_{\bf n})_{{\bf n}\in \Lambda}$ be the corresponding sequence of type II Hermite-Pad\'e approximants of ${\bf f}= \widehat{\bf s} + {\bf r}$. Then, for $j=1,\ldots m$
\begin{equation} \label{convunif} \lim_{{\bf n}\in \Lambda}
\frac{P_{{\bf n},j}}{Q_{{\bf n}}} = f_j = \widehat{s}_{1,j} + r_j, \qquad \mbox{inside} \qquad (\mathbb{C}\setminus \Delta_{1})^{\prime},
\end{equation}
the set obtained deleting from $\mathbb{C}\setminus \Delta_{1} $ the poles of all the $r_j$. For each $\varepsilon > 0$ sufficiently small, there exists $N > 0$ such that for all  ${\bf n} \in \Lambda$ with $\sum_{k=1}^m n_k \geq N$ we have $\deg Q_{\bf n} = \sum_{k=1}^m n_k$, if $\zeta$ is a pole of some $r_j$ of order $\kappa$ then $Q_{\bf n}$ has exactly $\kappa$ zeros in the disk $\{z : |z-\zeta| < \varepsilon\}$, and $Q_{\bf n}$ has exactly $\sum_{k=1}^m (n_k-d_k) $ simple zeros in $\stackrel{\circ}{\Delta}_1 $ (the interior of $\Delta_1$ with the Euclidean topology of $\mathbb{R}$).
\end{theorem}

The AT systems introduced in \cite{Nik} play an important role in the study of Nikishin systems. Given a multi-index ${\bf n}=(n_0,\ldots, n_m) \in \mathbb{Z}_+^{m+1} \setminus \{\bf 0\}$, a system of real continuous function $(u_0,\ldots,u_m)$ defined on an interval $\Delta$ is said to be an AT-system for ${\bf n}$ on $\Delta$ if for any  polynomials with real coefficients $p_0,\ldots,p_m$ such that  $\deg p_j \leq n_j -1,$ not all identically equal to zero, the linear form
\[
 {\mathcal{L}}(x)=p_0(x)u_0+\cdots+p_m(x) u_m(x)
\]
has at most $n_0+\cdots+n_m -1$ zeros on $\Delta$. If this is true for every $\bf n \in {\mathbb{Z}}_+^{m+1}\setminus \{\bf 0\}$ we have an AT system on $\Delta$.

According to \cite[Theorem 1.1]{KN}, if one fixes  $N < n_0+\cdots+n_m$ points $x_1,\ldots,x_N$ in the interior of $\Delta$ then we can find a convenient set of polynomials $p_0,\ldots,p_m$ with $\deg p_j \leq n_j -1, j=0,\ldots,m,$ such that $ {\mathcal{L}}$ changes sign at $x_1,\ldots,x_N$ and has no other sign change  in $\Delta.$

Nikishin systems are examples of AT systems, see \cite[Theorem 1.1]{FL4} and \cite[Theorem 1.1]{FL4II}. For the proof of Theorem \ref{TCTII} we need the AT property for polynomial modifications of Nikishin systems.

\begin{theorem}\label{teoAT}
Let $(s_{1,1},\ldots, s_{1,m})= \mathcal{N}(\sigma_{1},\ldots, \sigma_{m})$ and $\Lambda \subset \mathbb{Z}_+^{m+1}$ be given, where
\begin{equation}\label{cond3} \sup \left\{\max_{j=0,\ldots,m}(n_j) - \min_{k=0,\ldots,m}(n_k) : {\bf n} =(n_0,\ldots,n_m) \in \Lambda\right\} \leq C < \infty.
\end{equation}
Let $(t_{0},\ldots,t_{m})$ be a vector  polynomial with real coefficients whose zeros lie in $\mathbb{C}\setminus \Delta_{m}$ and for $j\neq k$ the zeros of $t_j$ and $t_k$ are distinct. Assume that $\mbox{\rm dist}(\Delta_{m-1},\Delta_m) > 0$  or $\sigma_m \in \mathcal{C}(\Delta_m)$. Then, on any interval $\Delta\subset \mathbb{R}\setminus \Delta_{1}$, $(t_0,t_1\hat{s}_{1,1}, \ldots, t_m\hat{s}_{1,m})$ is an AT system  for all ${\bf n} \in \Lambda$ with $n_0+\cdots+n_m$ sufficiently large.
\end{theorem}

The previous result is a consequence of the convergence of a multi-point version of type I Hermite Pad\'e approximants (see Theorem \ref{CTI} below).

In the sequel,   we maintain the notation introduced above.  The paper is organized as follows. Section 2 contains some auxiliary results and the proof of Theorem
\ref{TCTII} assuming that Theorem \ref{teoAT} holds. In Section 3 we study the convergence of multi-point type I Hermite Pad\'e approximants and prove Theorem \ref{teoAT}.

\section{Auxiliary results}

We need some preliminary facts about Nikishin systems which we will summarize below.
The following Lemma is easy to deduce. For a proof see \cite[Lemma 2.1]{LS}. In this section the  rational functions $r_j$ are allowed to have common poles except in the proof of Theorem \ref{TCTII}.

\begin{lemma}\label{psl}
Let ${\bf s} = (s_{1,1},\ldots,s_{1,m}) = \mathcal{N}(\sigma_1,\ldots,\sigma_m)$ be given. Assume that there exist polynomials with real coefficients $a_0,\ldots,a_m$ and a polynomial $w$ with real coefficients whose zeros lie in $\mathbb{C} \setminus \Delta_1$ such that
\[\frac{\mathcal{A}_0(z)}{w(z)} \in \mathcal{H}(\mathbb{C} \setminus \Delta_1)\qquad \mbox{and} \qquad \frac{\mathcal{A}_0(z)}{w(z)} = \mathcal{O}\left(\frac{1}{z^N}\right), \quad z \to \infty,
\]
where $\mathcal{A}_0  := a_0 + \sum_{k=1}^m a_k  \widehat{s}_{1,k} $ and $N \geq 1$. Let $\mathcal{A}_1  := a_1 + \sum_{k=2}^m a_k  \widehat{s}_{2,k} $. Then
\begin{equation} \label{eq:3}
\frac{\mathcal{A}_0(z)}{w(z)} = \int \frac{\mathcal{A}_1(x)}{(z-x)} \frac{d\sigma_1(x)}{w(x)}.
\end{equation}
If $N \geq 2$, we also have
\begin{equation} \label{eq:4}
\int x^{\nu}  \mathcal{A}_1(x)  \frac{d\sigma_1(x)}{w(x)} = 0, \qquad \nu = 0,\ldots, N -2.
\end{equation}
In particular, $\mathcal{A}_1$ has at least $N -1$ sign changes in  $\stackrel{\circ}{\Delta}_1 $.
\end{lemma}

In all that follows, when we write $\mathcal{O}(1/z^N)$ it is understood that $z \to \infty$ and the limit is taken  along any curve which is not tangent to the half straight line containing the support of the measures under consideration. By $T$ we denote the least common multiple of the denominators $t_1,\ldots,t_m$ of the rational functions $r_1,\ldots,r_m$. In this section ${\bf n} = (n_1,\ldots,n_m) \in \mathbb{Z}_+^m \setminus \{{\bf 0}\}$ and $|{\bf n}| := n_1 +\cdots+n_m$.

\begin{lemma} \label{OTII}
Let ${\bf R}_n=\left(P_{{\bf n},1}/Q_{{\bf n}},\ldots, P_{{\bf n},m}/Q_{{\bf n}}\right)$ be a type II Hermite-Pad\'e approximant with respect to ${\bf f} = \widehat{\bf s} + {\bf r}$ and ${\bf n} \in \mathbb{Z}_+^m \setminus \{{\bf 0}\}$. Assume that $n_j > D:= \deg T, j=1\ldots,m$. Then, for each $j=1,\ldots,m$
\begin{equation}\label{orto2}
\int  x^{\nu} t_j(x)Q_{{\bf n}}(x)ds_{1,j}(x)=0,  \qquad \nu=0,1, \ldots n_{j}-d_j-1,
\end{equation}
where $d_j = \deg t_j$. It follows that for any polynomials $p_j, \deg p_j \leq n_j - d_j -1, j=1,\ldots,m$
\begin{equation}\label{orto4}
\int  Q_{{\bf n}}(x) \left( p_1  t_1  + \sum_{j=2}^m p_j  t_j  \widehat{s}_{2,j} \right)(x) d\sigma_1(x)=0,
\end{equation}
and for any polynomials $p_j, \deg p_j \leq n_j - D -1, j=1,\ldots,m$
\begin{equation}\label{orto3}
\int  (TQ_{{\bf n}})(x) \left( p_1     + \sum_{j=2}^m p_j    \widehat{s}_{2,j} \right)(x) d\sigma_1(x)=0.
\end{equation}
Hence $ Q_{{\bf n}} $ has at least  $\left|{\bf n}\right|-mD$  zeros in   $\stackrel{\circ}{\Delta}_{1}$. Moreover,
\begin{equation}\label{restoint}
(T Q_{{\bf n}} \widehat{s}_{1,j} + TQ_{{\bf n}}(v_{j}/t_{j})  -T P_{{\bf n},j})(z)=\int \frac{(TQ_{{\bf n}})(x)}{z-x}ds_{1,j}(x).
\end{equation}
\end{lemma}

{\bf Proof.} Take $\bf n$ so that $n_j \geq D, j=1,\ldots,m$. For each $j=1,\ldots,m,$ ii) becomes
\[
Q_{{\bf n}}(z)\left(\hat{s}_{1,j}+\frac{v_{j}}{t_{j}}\right)(z) -P_{{\bf n},j}(z)=\mathcal{O}\left(\frac{1}{z^{n_{j}+1}}\right)\in {\mathcal{H}}\left({ {\mathbb{C}}}\setminus \left(\Delta_1 \cup \aleph_j\right)\right),
\]
where $\aleph_j$ is the set of zeros of $t_j$. Multiplying by $t_j$, we obtain
\[
\left(Q_{{\bf n}}t_j\hat{s}_{1,j}+Q_{{\bf n}}v_{j} -t_jP_{{\bf n},j}\right)(z)=\mathcal{O}\left(\frac{1}{z^{n_{j}-d_j+1}}\right)\in {\mathcal{H}}\left( {\mathbb{C}}\setminus \Delta_1 \right),
\]
and using  (\ref{eq:4}) in Lemma \ref{psl}, relation (\ref{orto2}) follows.
Taking  linear combinations of the relations given by (\ref{orto2}) we arrive at (\ref{orto4}).

In (\ref{orto2}) we can replace $x^{\nu}, \nu=0,\ldots,n_j -d_j -1,$ by $x^{\nu}T/t_j, \nu=0,\ldots,n_j -D -1$
and taking linear combinations of the orthogonality relations thus obtained  we get (\ref{orto3}).

From \cite[Theorem 1.1]{FL4II} we know that  $(1,\hat{s}_{2,2},\ldots, \hat{s}_{2,m})$ forms  an AT-system on $\Delta_{1}$. Suppose that $Q_{{\bf n}}$  has at most $\left|{\bf n}\right|-mD -1$ sign changes in   $\stackrel{\circ}{\Delta}_{1}$. Then, we can  choose $p_1.\ldots,p_m$ conveniently so that $p_1 + \sum_{j=2}^m p_j   \widehat{s}_{2,j}, \deg p_j \leq n_j - D -1, j=1,\ldots,m,$ has simple zeros at the points of sign change of $Q_{{\bf n}}$ on $\stackrel{\circ}{\Delta}_{1}$ and no other zero in $\Delta_1$. This contradicts (\ref{orto3}).  Thus, $Q_{{\bf n}}$  has at least  $\left|{\bf n}\right|-mD$ zeros  in   $\stackrel{\circ}{\Delta}_{1}$. (This type of argument will be used frequently and we will not detail it each time.) Finally, (\ref{restoint}) is a consequence of (\ref{eq:3}) in Lemma \ref{psl}.
\hfill $\Box$  \medskip

It is well known that given a measure $\sigma \in {\mathcal{M}}(\Delta),$ where $\Delta$ is contained in a half line, there exists a measure $\tau \in
{\mathcal{M}}(\Delta)$ and ${\ell}(z)=a z+b, a = 1/|\sigma|, b \in {\mathbb{R}},$ such that
\begin{equation} \label{s22}
{1}/{\widehat{\sigma}(z)}={\ell}(z)+ \widehat{\tau}(z),
\end{equation}
where $|\sigma|$ is the total variation of the measure $\sigma.$  See  \cite[Appendix]{KN} for measures with compact support and \cite[Lemma 2.3]{FL4} when the support is contained in a half line. $\tau$  is often called the inverse measure of $\sigma.$  Such measures appear frequently in our reasonings, hence we will fix a notation to distinguish them. In relation to measures  denoted by $s_{j,k}$ their inverse $\tau_{j,k}$ will carry over the corresponding sub-indices. The same goes for the  polynomials $\ell_{j,k}$. That is,
\[
{1}/{\widehat{s}_{j,k}(z)}  ={\ell}_{j,k}(z)+
\widehat{\tau}_{j,k}(z).
\]

For each $j \in \{2, \ldots ,m\}$ we define an auxiliary Nikishin
system
\[ S^j=(s_{2,2}^j, \ldots , s_{2,m}^j)= {\mathcal{N}}(\sigma_2^j,\ldots,\sigma_m^j) :=
\]
\[
{\mathcal{N}}(\tau_{2,j},\widehat{s}_{2,j} d \tau_{3,j}, \ldots ,
\widehat{s}_{j-1,j} d \tau_{j,j}, \widehat{s}_{j,j} d
\sigma_{j+1}, \sigma_{j+2}, \ldots , \sigma_m)\,.
\]
We also define $S^1=(s_{2,2}^1, \ldots , s_{2,m}^1) = {\mathcal{N}}(\sigma_2,\ldots,\sigma_m)$.

These auxiliary systems were used in the proof of \cite[Lemmas 5-6]{Bus} (see also \cite[Theorem 3.1.3]{DrSt2}).
Subsequently,
in \cite[Lemma 3.2]{FL4}, several formulas involving ratios of Cauchy transforms were proved (when the supports of the generating measures are unbounded or overlapping see \cite[Lemma 2.10]{FL4II} instead).
For convenience of the reader, we
write those formulas using the notation introduced above. We have:

\begin{equation}\label{eq:k}
\frac{1}{\widehat{s}_{2,j}} = \ell_{2,j} +
\widehat{\tau}_{2,j}\,,
\end{equation}

\begin{equation} \label{4.4}
\frac{\widehat{s}_{1,k}}{\widehat{s}_{1,1}} =
\frac{|s_{1,k}|}{|s_{1,1}|} - \langle \tau_{1,1},\langle s_{2,k},\sigma_1
\rangle \widehat{\rangle}  , \qquad  1 < k \leq m,
\end{equation}

\begin{equation}\label{eq:l}
\frac{\widehat{s}_{2,k} }{\widehat{s}_{2,j} } = b_{j,k} +
  (-1)^{k-1}\widehat{s}_{2,k+1}^j  + c_{j,k} \widehat{s}_{2,k}^j  \,, \qquad k=
2,\ldots,j-1\,,
\end{equation}
and
\begin{equation}\label{eq:m}
\frac{\widehat{s}_{2,k}}{\widehat{s}_{2,j}} = b_{j,k} +
c_{j,k}\widehat{s}_{2,k}^j(z) \,, \qquad k =j+1,\ldots,m\,,
\end{equation}
where the $b_{j,k} $ and $c_{j,k}$ denote (perfectly determined) constants.

\begin{definition} \label{asociado} Let ${\bf n}=(n_1, \ldots, n_m) \in
{\mathbb{Z}}_+^m$. For each $j=1, \ldots, m,$ we define an
associated multi-index ${{\bf n}}^j=(n^j_{2}, \ldots, n^j_{m})$
whose $m-1$ components are
\[
n^j_k= \left\{ \begin{array}{lcl} \min (n_1, \ldots , n_{k-1}, n_j-1),& \mbox{when } & k=2, \ldots ,j,\\ & \\
                                     \min (n_j,  n_k), & \mbox{when } & k=j+1, \ldots, m\,.
\end{array}
\right.
\]
We denote $|{{\bf n}}^j|= \sum_{k=2}^m n^j_{k} $.
\end{definition}

For $j=1,,\ldots,m$, set $\Phi_{{\bf n},j}:=(T Q_{{\bf n}} \widehat{s}_{1,j} + TQ_{{\bf n}}(v_{j}/t_{j})  -T P_{{\bf n},j}) $ (see (\ref{restoint})).
\begin{lemma} \label{lem:b}
Let ${{\bf n}} = (n_1,\ldots,n_m)$ be a multi-index such that $n_k > D = \deg T, k=1,\ldots,m$. Then, for each $j=1,
\ldots, m$
\begin{equation}\label{eq:j}
\int  (p_k  \Phi_{{{\bf n}},j} )(x) ds_{2,k}^j(x)  = 0,  \quad \deg p_k < n_{k}^{j}-D, \quad k=2,\ldots,m.
\end{equation}
Set $N^{j}=\left|{\bf n}^{j}\right|+n_{j}$. For each $j=1,\ldots,m$, there exists a monic polynomial $w_{{\bf n},j}$, $\deg w_{{\bf n},j}\geq  \left|{\bf n}^{j}\right|-(m-1)D$, with simple zeros which lie in  $\stackrel{\circ}{\Delta}_{2}$, such that
\begin{equation}\label{eq:p:1}
\frac{\Phi_{{{\bf n}},j}(z)}{w_{{\bf n},j}(z)}
=\mathcal{O}\left(\frac{1}{z^{N^{j}-mD+1}}\right) \in H(\mathbb{C}\setminus \Delta_1)\,,
\end{equation}
and
\begin{equation}\label{ROIK}
0=\int x^{\nu}T(x)Q_{{\bf n}}(x)\frac{ds_{1,j}(x)}{w_{{\bf n},j}(x)}, \qquad \nu=0,1,\ldots N^{j}-mD-1.
\end{equation}
\end{lemma}

{\bf Proof} Fix $j \in \{1,\ldots,m\}$. From the definition of
$\Phi_{{\bf n},j} $ and (\ref{restoint}),  we obtain
\[
\int  (p_k  \Phi_{{{\bf n}},j} )(x) ds_{2,k}^j(x)=\int
p_k(x) \int  \frac{T(t)Q_{{\bf n}}(t)}{x-t} d s_{1,j}(t) ds_{2,k}^j(x).
\]
Since $\deg p_k < n_k^j -D\leq n_j-D,$ from (\ref{orto2}) and
Fubini's theorem, it follows that
\[
\int  p_k(x) \int  \frac{T(t)Q_{{\bf n}}(t)}{x-t} d
s_{1,j}(t) ds_{2,k}^j(x)=
\]
\begin{equation}
\label{fundemo}\int  \int
\frac{(p_kTQ_{{\bf n}})(t)}{x-t} d s_{1,j}(t)ds_{2,k}^j(x)= -\int (p_kTQ_{{\bf n}}
\widehat{s}_{2,k}^j)(t) d s_{1,j}(t)\,.
\end{equation}

First, we prove the statement of the lemma for the case $j+1 \leq
k \leq m$. If $j = m$, the set of $k$ is empty and there is nothing to prove. Let $j \leq m-1$.
Using (\ref{eq:m}), we obtain
\[
- \int  (p_k TQ_{{\bf n}} \widehat{s}_{2,k}^j)(t)d
s_{1,j}(t) =
  \int (p_k TQ_{{\bf n}})(t)
\displaystyle{\left(\frac{\widehat{s}_{2,k}(t)}{\widehat{s}_{2,j}(t)}
- a_{j,k} \right)} d s_{1,j}(t) =
\]
\[-a_{j,k} \int  (p_kTQ_{{\bf n}})(t) d s_{1,j}(t)+ \int
 (p_kTQ_{{\bf n}})(t) d s_{1,k}(t)\,.
\]
By hypothesis
$\deg p_k < n_k^j-D \leq \min \{n_j-D,n_k-D\}$. Taking (\ref{orto2}) into account, we deduce that both terms after the last equality vanish. Hence the
first case is proved.

Now, we analyze the case when $2 \leq k \leq j$. Using several
times formula (\ref{eq:l}) to make $k$ descend to $2$ and finally
formula (\ref{eq:k}), we obtain the equalities
\begin{equation} \label{eq:tio6} \widehat{s}^j_{2,k}  =
(-1)^k\left(\frac{\widehat{s} _{2,k-1}}{\widehat{s}_{2,j}}- a_{j,k-1} - c_{j,k-1}\widehat{s}^j_{2,k-1}\right) =
\end{equation}
\[ (-1)^k\left(\frac{\widehat{s} _{2,k-1}}{\widehat{s}_{2,j}}- a_{j,k-1}- (-1)^{k-1}c_{j,k-1}\left( \frac{\widehat{s} _{2,k-2} }{\widehat{s}_{2,j}}- a_{j,k-2} - c_{j,k-2}\widehat{s}^j_{2,k-2}\right)
\right) =\]
\[ \cdots =
{\mathcal{L}}_j^* + \frac{1}{\widehat{s}_{2,j}}\sum_{l=1}^{k-1} c_{l}^*
 \widehat{s}_{2,l}\,,
\]
where ${\mathcal{L}}_{j}^*$ denotes a polynomial of degree $1$,
$\widehat{s}_{2,1} \equiv 1$, and the $c_{l}^*, l= 1,\ldots,k-1,$ are
constants.
Sustituting (\ref{eq:tio6}) into (\ref{fundemo}), we obtain
\[
\int  (p_kTQ_{{\bf n}}\widehat{s}_{2,k}^j)(t) d
s_{1,j}(t) =
\]
\[
 -\int (p_kTQ_{{\bf n}})(t) {\mathcal{L}}_{j}^*(t) d
s_{1,j}(t) - \sum_{l=1}^{k-1}c_{i}^* \int  (p_kTQ_{{\bf n}})(t)d
s_{1,l}(t)\,.
\]
From hypothesis $\deg p_k   \leq
\min(n_1-D-1,\ldots,n_{k-1}-D-1,n_j-D-2)$ and using (\ref{orto2})  it follows that all the
integrals on the right hand side of this equality are zero. Hence, (\ref{eq:j}) holds.

From  (\ref{eq:j}) it follows that for any polynomials $p_k, \deg p_k < n_k^j-D, k=2, \ldots ,m$
\begin{equation} \label{eq:i:2}
\int   \Phi_{{\bf n},j}(x)\left(p_2 + \sum_{k=3}^m
p_k \widehat{s}_{3,k}^j\right)(x) d \sigma_2^j(x) =0.
\end{equation}
According to  \cite[Theorem 1.1]{FL4II}, $(1, \widehat{s}_{3,3}^j, \ldots ,\widehat{s}_{3,m}^j)$
forms an AT system on $\Delta_2$. Thus, using (\ref{eq:i:2}) it follows that
$\Phi_{{{\bf n}},j} $ has at least $|{{\bf n}}^j|-(m-1)D$ sign changes in
$\stackrel{\circ}{\Delta}_2$. Let $w_{{\bf n},j}$ be the monic
polynomial whose zeros are the points where $\Phi_{{{\bf n}},j} $
changes sign in the interior of $\Delta_2$.
Taking into account that $\deg w_{{\bf n},j}
\geq |{{\bf n}}^j|-(m-1)D$, we obtain (\ref{eq:p:1}) which together with
(\ref{eq:4}) implies (\ref{ROIK}).
\hfill $\Box$
\medskip

We are ready to prove convergence of type II Hermite Pad\'e approximants in a weaker sense. We use the concept of convergence in Hausdorff content.

Let $B$ be a subset of the complex plane $\mathbb{C}$. By
$\mathcal{U}(B)$ we denote the class of all coverings of $B$ by at
most a numerable set of disks. Set
$$
h(B)=\inf\left\{\sum_{j=1}^\infty
|U_j|\,:\,\{U_j\}\in\mathcal{U}(B)\right\},
$$
where $|U_j|$ stands for the radius of the disk $U_j$. The quantity
$h(B)$ is called the $1$-dimensional Hausdorff content of the
set $B$.

Let $(\varphi_n)_{n\in\mathbb{N}}$ be a sequence of complex functions
defined on a domain $\Omega \subset\mathbb{C}$ and $\varphi$ another
function defined on $\Omega$ (the value $\infty$ is permitted). We say that
$(\varphi_n)_{n\in\mathbb{N}}$ converges in Hausdorff content to
$\varphi$ inside $\Omega$ if for each compact
subset $\mathcal{K}$ of $\Omega$ and for every $\varepsilon
>0$, we have
$$
\lim_{n\to\infty} h\{z\in K :
|\varphi_n(z)-\varphi(z)|>\varepsilon\}=0
$$
(by convention $\infty \pm \infty = \infty$). We denote this writing $h$-$\lim_{n\to \infty} \varphi_n =
\varphi$ inside $\Omega$. Under certain assumptions, there is a lemma of A.A. Gonchar \cite[Lemma 1]{Gon} which allows to derive uniform convergence from convergence in Hausdorff content. We will follow this approach.

Let $s \in \mathcal{M}(\Delta)$ where $\Delta$ is contained in a half line of the real axis and let $r=p/q, (p,q) \equiv 1, \deg p < \deg q,$ be a rational function with real coefficients whose poles belong to $\mathbb{C} \setminus \Delta$. Set $f = \widehat{s} +r$. Fix $\kappa \geq -1$ and $\ell \in \mathbb{Z}_+$. Consider a sequence of polynomials $(w_n)_{n \in \Lambda},  \Lambda \subset \mathbb{Z}_+,$ such that $\deg w_n = \kappa_n \leq 2n + \kappa - \ell +1$, whose zeros lie in $\mathbb{R} \setminus \Delta$. Let $(R_n)_{n \in \Lambda}$ be a sequence of rational functions $R_n = p_n/q_n$ with real coefficients satisfying the following conditions for each $n \in \Lambda$:
\begin{itemize}
\item[a)] $\deg p_n \leq n + \kappa,\quad  \deg q_n \leq n, \quad q_n \not\equiv 0,$
\item[b)] $ {(q_n f - p_n)(z)}/{w_n}(z) = \mathcal{O}\left( {z^{-n-1 + \ell}}\right) \in \mathcal{H}(\mathbb{C}\setminus (\Delta \cup [r=\infty]) $.
\end{itemize}
We call $R_n$ an incomplete diagonal multi-point Pad\'e approximant of $f$.

Notice that in this construction for each $n \in \Lambda$ the number of free parameters equals $2n + \kappa +2$ whereas the number of homogeneous linear equations to be solved  to find $q_n$ and $p_n$ equals $2n + \kappa - \ell + 1$. When $\ell =0$ there is only one more parameter than equations, $R_n$ is defined uniquely, and coincides with a diagonal multi-point Pad\'e approximant of $f$. When $\ell \geq 1$ uniqueness is not guaranteed, thus the term incomplete.

For sequences of incomplete diagonal multi-point Pad\'e approximants, the following Stieltjes type theorem was proved in \cite[Lemma 2]{Bus} in terms of convergence in logarithmic capacity. The proof using Hausdorff content is basically the same.

\begin{lemma} \label{BusLop}
Assume that $(R_n)_{n \in \Lambda}$ satisfies {\rm a)-b)} and either the number of zeros of $w_n$ lying on a bounded segment of $\mathbb{R} \setminus \Delta$ tends to infinity as $n\to\infty, n \in \Lambda$, or $s \in \mathcal{C}(\Delta)$.
Then
\[ h-\lim_{n \in \Lambda} R_n = \widehat{s},\qquad \mbox{inside}\qquad  \mathbb{C} \setminus \Delta.
\]
\end{lemma}
\medskip

An immediate consequence of this lemma is the following.

\begin{lemma} \label{convHaus2} Let $ (s_{1,1},\ldots,s_{1,m}) = \mathcal{N}(\sigma_1,\ldots,\sigma_m)$ and $  (r_1,\ldots,r_m)$ be given,  where the rational functions $r_j, j=1,\ldots,m$ have real coefficients and  their poles  lie in  ${\mathbb{C}}\setminus \Delta_1$. Let $\Lambda \subset \mathbb{Z}_+^{m}$ be an infinite sequence of distinct multi-indices  satisfying $(\ref{cond1})$.
Assume that either $\mbox{\rm dist}(\Delta_{1},\Delta_2) >0$  or $\sigma_1\in \mathcal{C}(\Delta_1)$. Then, for $j=1,\ldots,m$
\[ h-\lim_{{\bf n}\in \Lambda}
\frac{P_{{\bf n},j}}{Q_{{\bf n}}} = f_j = \widehat{s}_{1,j} + r_j, \qquad \mbox{inside} \qquad \mathbb{C}\setminus \Delta_{1}.
\]
\end{lemma}

{\bf Proof.} Fix $j =1,\ldots,m$. According to Lemma \ref{lem:b}, the zeros of $w_{{\bf n},j}$ lie in $\mathbb{R} \setminus \Delta_1$. Using (\ref{eq:p:1}) and (\ref{cond1}), it is easy to verify that  for all sufficiently large the rational function $P_{{\bf n},j}/Q_{\bf n}$ is an incomplete Pad\'e approximant of $f_j$ taking $n = |\bf n|, \kappa = -1$ and choosing $\ell$ appropriately. It is easy to check as well that  $\sigma_1 \in \mathcal{C}(\Delta_1)$ implies $s_{1,j} \in \mathcal{C}(\Delta_1)$. Therefore, the statement readily follows from Lemma \ref{BusLop}. \hfill $\Box$
\medskip

{\bf Proof of Theorem \ref{TCTII}.}
We assume that Theorem \ref{teoAT} holds. Therefore, there exists an $N$ such that $(t_{1},t_{2}\hat{s}_{2,2},\ldots,t_{m}\hat{s}_{2,m})$ forms an AT-system on $\Delta_{1}$ for the multi-index $(n_{1}-d_{1},n_{2}-d_{2},\ldots,n_m-d_m)$ if $ \left|{\bf n}\right|>N$. From (\ref{orto4}) it follows that
$Q_{{\bf n}}$ has at least  $|{\bf n}|- D, D := d_{1}+d_{2}+\ldots+d_m$ zeros in $\stackrel{\circ}{\Delta}_{1}$  if $ \left|{\bf n}\right|>N$.

Let $\zeta$ be a pole of $r_j$ of order $\kappa$. Using Lemma \ref{convHaus2} and \cite[Lemma 1]{Gon}) it follows that for any $\varepsilon > 0$ sufficiently small $Q_{\bf n}$ has at least $\kappa$ zeros in $\{z: |z-\zeta| < \varepsilon\}$. The poles of the different $r_j$ are distinct. Thus, for each $\varepsilon > 0$ and all sufficiently large $|{\bf n}|$ the zeros of $Q_{\bf n}$ are either in $\stackrel{\circ}{\Delta}_{1}$ or inside an $\varepsilon$ neighborhood of the poles of the $r_j$. Lemma \ref{convHaus2} and \cite[Lemma 1]{Gon}) then imply (\ref{convunif}).  \hfill $\Box$

\section{Type I Hermite-Pad\'e approximation}\label{CTI}

In this section ${\bf n} := (n_0,\ldots,n_m) \in \mathbb{Z}_+^{m+1} \setminus \{{\bf 0}\}, |{\bf n}| = n_0+\cdots+n_m,$ and $N_{\bf n} = \max\{n_0,n_1-1,\ldots,n_m-1\}$.

We need an analogue of Lemma \ref{convHaus2} for type I Hermite-Pad\'e approximation. Let $(s_{1,1},\ldots,s_{1,m}) = \mathcal{N}(\sigma_1,\ldots,\sigma_m)$, and ${\bf n} \in \mathbb{Z}_+^{m+1}  \setminus \{\bf 0\}$ be given.  Fix   $\ell \in \mathbb{Z}_+$ and  a polynomial $w_{\bf n}, \deg w_{\bf n} \leq |{\bf n}| -\ell-1,$ whose zeros lie in $\mathbb{R} \setminus\Delta_1$.  Consider a vector polynomial $\left(p_{{\bf n},0},\ldots, p_{{\bf n},m}\right),$ not  identically equal zero, which satisfies:
\begin{itemize}
\item[i')] $\deg p_{{\bf n},j} \leq n_j -1, j=0,\ldots,m,$
\item[ii')] $\mathcal{A}_{{\bf n},0}/w_{\bf n} = \mathcal{O}(1/z^{|{\bf n}|-N_{\bf n} -\ell}) \in \mathcal{H}(\mathbb{C} \setminus \Delta_1), \quad
    \mathcal{A}_{{\bf n},0} := p_{{\bf n},0} + \sum_{k= 1}^m p_{{\bf n},k} \widehat{s}_{1,k}.$

\end{itemize}
We call  $\left(p_{{\bf n},0},\ldots, p_{{\bf n},m}\right)$ an incomplete type I multi-point Hermite-Pad\'e  approximation of $(\widehat{s}_{1,1},\ldots,\widehat{s}_{1,m})$ with respect to $w_{\bf n}$ and $\bf n$.

The proof of the following lemma reduces to Lemma \ref{BusLop}. When $\ell = 0$ it is contained in \cite[Lemma 3.1]{LS} and for general $\ell$ it requires no substantial change so we omit it.

\begin{lemma} \label{CCTI}  Let ${\bf s}= (s_{1,1},\ldots,s_{1,m}) = \mathcal{N}(\sigma_1,\ldots,\sigma_m)$ and $\Lambda \subset \mathbb{Z}_+^{m}$ be an infinite sequence of distinct multi-indices. Fix  $\ell \in \mathbb{Z}_+$ and a sequence $(w_{\bf n})_{{\bf n}\in \Lambda}, \deg w_{\bf n} \leq |{\bf n}|  - \ell -1,$  of polynomials  whose zeros lie in $\mathbb{R} \setminus\Delta_1$.  Consider a sequence of incomplete type I multi-point Hermite-Pad\'e approximants of ${\bf s}$ with respect to $(w_{\bf n})_{{\bf n} \in \Lambda}$.
Assume that $(\ref{cond3})$ takes place and that either $\mbox{\rm dist}(\Delta_{m-1},\Delta_m) >0$ or $\sigma_m \in \mathcal{C}(\Delta_m)$. Then,
for each fixed $j=0,\ldots,m-1$
\begin{equation} \label{convHaus}
h-\lim_{{\bf n}\in \Lambda}\frac{p_{{\bf n}, j}}{p_{{\bf n},m}} = (-1)^{m-j}\widehat{s}_{m,j+1}, \quad h-\lim_{n\in \Lambda}\frac{p_{{\bf n}, m}}{p_{{\bf n},j}} = \frac{(-1)^{m-j}}{\widehat{s}_{m,j+1} },
\end{equation}
inside $\mathbb{C} \setminus \Delta_m$.  There exists a constant $C_1$, independent of $\Lambda$, such that  for all ${\bf n} \in \Lambda,$ the polynomials $p_{{\bf n},j}, j=0,\ldots,m,$  have at least $(|{\bf n}|/(m+1)) - C_1$ zeros in $\stackrel{\circ}{\Delta}_m $.
\end{lemma}

We are interested in a special type of incomplete type I multi-point Hermite Pad\'e approximant.
Let $(t_0,t_{1},\ldots,t_{m}), \deg t_j = d_j,$ be a vector polynomial with real coefficients  and let $D:=\sum_{j=0}^{m}d_{j}$.
Fix  ${{\bf n}} = (n_0,\ldots,n_m) \in \mathbb{Z}_+^{m+1}, n_{j}>D, j=0,\ldots,m$. Given   $(s_{1,1},\ldots, s_{1,m})= \mathcal{N}(\sigma_{1},\dots, \sigma_{m})$ and a
polynomial with real coefficients $w_{{\bf n}}, \deg w_{{\bf n}} \leq |{{\bf n}}|-D-1,$  whose zeros lie in  $\mathbb{R}\setminus \Delta_{1}$, there exist polynomials $a_{{{\bf n}},0},a_{{{\bf n}},1}, \ldots, a_{{{\bf n}},m}$, not all identically equal to zero, such that:
\begin{itemize}
\item[a')] $\deg a_{{{\bf n}},j}t_{j} \leq n_j-1, j=0,\ldots,m,$
\item[b')] $\frac{a_{{{\bf n}},0}(z)t_0(z)+\sum_{j=1}^m a_{{{\bf n}},j}(z)  t_j(z)\hat{s}_{1,j}(z) }{w_{{\bf n}}(z)} = {\mathcal{O}}(1/z^{|{{\bf n}}|- N_{\bf n}-D}) \in \mathcal{H}(\mathbb{C}\setminus \Delta_{1})$.
\end{itemize}

In other words, we use the freedom in the construction of the incomplete type I Hermite-Pad\'e approximants to force the polynomials $p_{{\bf n},j}, j=0,\ldots,m$ to have some predetermined zeros. $(a_{{\bf n},0},\ldots,a_{{\bf n},m})$ could be regarded as a type I multi-point Hermite Pad\'e approximant of $(t_1\widehat{s}_{1,1}/t_0,\ldots,t_m\widehat{s}_{1,m}/t_0)$ with respect to $(n_0-d_0,\ldots,n_m-d_m)$.

\begin{theorem}\label{CTI}
Let $(s_{1,1},\ldots,s_{1,m}) = \mathcal{N}(\sigma_1,\ldots,\sigma_m)$ and $\Lambda \subset \mathbb{Z}_+^{m+1}$ be an infinite sequence of distinct multi-indices which verifies $(\ref{cond3})$.  Assume  that the polynomials $t_{0},t_{1},\ldots t_{m}$  have no common zeros, they all  lie in $\mathbb{C}\setminus \Delta_{m},$
and either $\mbox{\rm dist}(\Delta_{m-1},\Delta_m)> 0$  or $\sigma_m \in \mathcal{C}(\Delta_m)$. Suppose that for each ${\bf n} \in \Lambda$ the polynomials $a_{{{\bf n}},0},a_{{{\bf n}},1}, \ldots, a_{{{\bf n}},m}$ satisfy {\rm a')-b')}.
Then, for $j=0,1,\ldots m-1$
\begin{equation} \label{cu} \displaystyle \lim_{{\bf n}\in \Lambda}\frac{a_{{{\bf n}},j}}{a_{{{\bf n}},m}}= (-1)^{m-j}\frac{t_m}{t_j}\hat{s}_{m,j+1},\qquad  \mbox{inside} \qquad (\mathbb{C}\setminus \Delta_{m})^{\prime},
\end{equation}
the set obtained deleting from $\mathbb{C}\setminus \Delta_{m} $ the zeros of all the polynomials $t_j.$ There exists a constant $C_1$, independent of $\Lambda$, such that for each $j=0,\ldots,m$ and ${\bf n} \in \Lambda,$ the polynomials $a_{{\bf n},j}$  have at least $(|{\bf n}|/(m+1)) - C_1$ zeros in $\stackrel{\circ}{\Delta}_m $. Fix $j,k=0,\ldots,m$. Let $\zeta$ be a zero of $t_k, k\neq j,$ of multiplicity $\kappa$. Then, for each $\varepsilon >0$ sufficiently small there exists an $N$ such that for all ${\bf n} \in \Lambda, |{\bf n}| > N$, $a_{{\bf n},j}$ has exactly $\kappa$ zeros in $\{z:|z-\zeta| < \varepsilon\}$. The remaining zeros of $a_{{\bf n},j}$ either lie on $\Delta_m$ or accumulate on $\Delta_m \cup \{\infty\}$ as $|{\bf n}| \to \infty$.
\end{theorem}

{\bf Proof.} The lower bound on the number of zeros of the $a_{{\bf n},j}$ in $\stackrel{\circ}{\Delta}_m $ is taken from Lemma \ref{CCTI}. Let $\overline{\jmath}$ be the last component of $(n_0,\ldots,n_m)$ such that $n_{\overline{\jmath}} = \min_{j=0,\ldots,m} (n_j)$.  Let us prove that  $a_{{\bf n},\overline{\jmath}}$ has at least $n_{\overline{\jmath}}-D -1$  zeros in $\stackrel{\circ}{\Delta}_{m} $. That is, for this component we want a more precise lower bound on the number of zeros.

From  \cite[Theorem 1.3]{FL4} (see also \cite[Theorem 3.2]{FL4II}), we know that  there exists a permutation $\lambda$ of $(0,\ldots,m)$ which reorders the components of $(n_0,n_1,\ldots,n_m)$ decreasingly, $n_{\lambda(0)} \geq \cdots \geq n_{\lambda(m)},$ and an associated Nikishin system $(r_{1,1},\ldots,r_{1,m}) = {\mathcal{N}}(\rho_{1},\ldots,\rho_m)$ such that
\[  \mathcal{A}_{{\bf n},0} = (q_{{\bf n},0} + \sum_{k=1}^m q_{{\bf n},k} \widehat{r}_{1,k})\widehat{s}_{1,\lambda(0)}, \quad \deg q_{{\bf n},k} \leq n_{\lambda(k)} -1, \quad k=0,\ldots,m,
\]
where $\mathcal{A}_{{\bf n},0} = a_{{{\bf n}},0} t_0 +\sum_{j=1}^m a_{{{\bf n}},j}   t_j \hat{s}_{1,j}  $ and $\widehat{s}_{1,\lambda(0)} \equiv 1$ when $\lambda(0) = 0$.
The permutation may be taken so that for all $0 \leq j <k \leq n$ with $n_j = n_k$ then also $\lambda(j) < \lambda(k)$. In this case,  see formulas (31) in the proof of \cite[Lemma 2.3]{FL4}, it follows that $q_{{\bf n},m}$ is either $a_{{\bf n},\overline{\jmath}}t_{\overline{\jmath}}$ or $-a_{{\bf n},\overline{\jmath}}t_{\overline{\jmath}}$.

Set
\[\mathcal{Q}_{{\bf n},j} := q_{{\bf n},j} + \sum_{k=j+1}^m q_{{\bf n},k} \widehat{r}_{1,k}, \quad j=0,\ldots,m-1, \quad \mathcal{Q}_{{\bf n},m} := q_{{\bf n},m}.\]
Suppose that $\lambda(0) = 0$ and thus $\widehat{s}_{1,\lambda(0)} \equiv 1$. Then, $N_{\bf n} = n_0 = n_{\lambda(0)}$. Due to b'), it follows that
\begin{equation} \label{resto} \frac{\mathcal{Q}_{{\bf n},0}(z)}{w_{\bf n}(z)} = \mathcal{O}(1/z^{|{\bf n}| - n_{\lambda(0)} -D}) \in \mathcal{H}(\mathbb{C}\setminus \Delta_1).
\end{equation}
When $\lambda(0) \neq 0$ we have that $N_{\bf n} = n_{\lambda(0)} - 1$ and $\widehat{s}_{1,\lambda(0)} = \mathcal{O}(1/z)$. Therefore, from b') we again have (\ref{resto}). Using (\ref{eq:4}), it follows that
\[ \int x^{\nu}  \mathcal{Q}_{{\bf n},1}(x)  \frac{d\rho_1(x)}{w_{\bf n}(x)} = 0, \qquad \nu = 0,\ldots, |{\bf n}| - n_{\lambda(0)}- D -2.
\]
This implies that $\mathcal{Q}_{{\bf n},1}$ has at least $|{\bf n}|  - n_{\lambda(0)}- D -1$ sign changes in $\stackrel{\circ}{\Delta}_1$. Let $w_{{\bf n},1}$ be the monic polynomial whose zeros are the points where $\mathcal{Q}_{{\bf n},1}$ changes sign in $\stackrel{\circ}{\Delta}_1$. Then
\[ \frac{\mathcal{Q}_{{\bf n},1}(z)}{w_{{\bf n},1}(z)} = \mathcal{O}\left(1/z^{|{\bf n}|  - n_{\lambda(0)} -n_{\lambda(1)} - D}\right) \in \mathcal{H}(\mathbb{C} \setminus \Delta_2).
\]
Using again (\ref{eq:4}) we get
\[ \int x^{\nu}  \mathcal{Q}_{{\bf n},2}(x)  \frac{d\rho_2(x)}{w_{{\bf n},1}(x)} = 0, \quad \nu = 0,\ldots, |{\bf n}|  - n_{\lambda(0)}- n_{\lambda(1)} - D -2,
\]
which implies that $\mathcal{Q}_{{\bf n},2}$ has at least $|{\bf n}|  - n_{\lambda(0)}- n_{\lambda(1)} - D -1$ sign changes in $\stackrel{\circ}{\Delta}_2$. Repeating the arguments $m$ times, it follows that $\mathcal{Q}_{{\bf n},m} = q_{{\bf n},m}$ has at least $n_{\overline{\jmath}} - D -1$ sign changes in $\stackrel{\circ}{\Delta}_m$ which implies the statement because $q_{{\bf n},m} = \pm a_{{\bf n},\overline{\jmath}}t_{\overline{\jmath}}$ and the zeros of $t_{\overline{\jmath}}$ are outside $\Delta_m$.

By Lemma \ref{CCTI}, for $j=0,\ldots m-1$,  we have
\begin{equation}\label{CC1}
h-\lim_{{\bf n}\in \Lambda}\frac{a_{{\bf n}, j} }{a_{{\bf n},m} } = (-1)^{m-j}\frac{t_m\widehat{s}_{m,j+1}}{t_j}, \qquad \mbox{inside} \qquad \mathbb{C} \setminus \Delta_m.
\end{equation}
From (\ref{CC1}) and \cite[Lemma 1]{Gon} it follows that each zero $\zeta$ of $t_j $ of multiplicity $\kappa$ attracts at least $\kappa$ zeros of $a_{{\bf n},m}$ when $|{\bf n}| \to \infty, {\bf n} \in \Lambda$ (recall that $t_m$ and $t_j$ are relatively prime). Let us show that the number of zeros attracted by $\zeta$ equals $\kappa$ and the rest of the zeros of $a_{{\bf n},m}$ either lie in $\stackrel{\circ}{\Delta}_m$ or accumulate on $\Delta_m \cup \{\infty\}$. Then \cite[Lemma 1]{Gon} and (\ref{CC1}) imply (\ref{cu}).

The index $\overline{\jmath}$ as defined above may depend on  ${\bf n} \in \Lambda$. Given $\overline{\jmath} \in \{0,\ldots,m\}$, let $\Lambda(\overline{\jmath})$ denote the set of all ${\bf n} \in \Lambda$ such that $\overline{\jmath}$ is the last component of $(n_0,\ldots,n_m)$ satisfying $n_{\overline{\jmath}} = \min_{j=0,\ldots,m} (n_j)$. Fix $\overline{\jmath}$ and suppose that $\Lambda(\overline{\jmath})$ contains infinitely many multi-indices.
Should $\overline{\jmath} = m$, then $a_{{\bf n},m}$ has $n_m -D -1$ zeros in $\stackrel{\circ}{\Delta}_m$ and the rest of its zeros converge to the zeros of the $t_j, j=0,\ldots,m-1,$ according to their multiplicity as needed.

Now, consider that  $ \overline{\jmath} \neq m$. By Lemma \ref{CCTI} we have
\begin{equation}\label{CC2}
h-\lim_{{\bf n}\in \Lambda(\overline{\jmath})}\frac{a_{{\bf n}, m}}{a_{{\bf n},\overline{\jmath}}} = (-1)^{m-\overline{\jmath}}\frac{t_{\overline{\jmath}}}{t_{m}\widehat{s}_{m,\overline{\jmath}+1}}\qquad \mbox{inside} \qquad \mathbb{C} \setminus \Delta_m.
\end{equation}
For each $j \in \{0,\ldots,m-1\} \setminus \{\overline{\jmath}\}$ and each zero $\zeta$ of multiplicity $\kappa$ of $t_j$ choose $\kappa$ zeros of $a_{{\bf n},m}$ that converge to $\zeta$ as $|{\bf n}| \to \infty, {\bf n} \in \Lambda(\overline{\jmath})$. Let $q_n$ be a monic polynomial with this set of points as its zeros. Obviousty, $\lim_{{\bf n} \in \Lambda(\overline{\jmath})} q_{\bf n}= \prod_{k=0}^{m-1}t_k/t_{\overline{\jmath}}$ (uniformly on compact subsets). From (\ref{CC2}), we get
\begin{equation}\label{CC5}
h-\lim_{n\in \Lambda(\overline{\jmath})}\frac{a_{{\bf n}, m}}{q_{\bf n}a_{{\bf n},\overline{\jmath}}} = (-1)^{m-\overline{\jmath}}\frac{t^2_{\overline{\jmath}}}{\prod_{k=0}^m t_{k}\widehat{s}_{m,\overline{\jmath}+1}}\qquad \mbox{inside} \qquad \mathbb{C} \setminus \Delta_m.
\end{equation}
Applying once more \cite[Lemma 1]{Gon}, it follows that for $j\in \{0,\ldots,m\} \setminus \{{\overline{\jmath}}\}$ each zero of $t_j$ attracts exactly as many zeros of $a_{{\bf n},\overline{\jmath}}$ as its multiplicity (notice that the zeros of $q_{\bf n}$ are canceled by zeros of $a_{{\bf n},m}$). Therefore, all the zeros of $a_{{\bf n},\overline{\jmath}}$ are located either on $\Delta_m$ or on a sufficiently small neighborhood of the zeros of the polynomial $\prod_{k=0}^m t_{k}/t_{\overline{\jmath}}$. Consequently,
\begin{equation}\label{CC6}
\lim_{{\bf n}\in \Lambda(\overline{\jmath})}\frac{a_{{\bf n}, m}}{a_{{\bf n},\overline{\jmath}}} = (-1)^{m-\overline{\jmath}}\frac{t_{\overline{\jmath}}}{t_{m}\widehat{s}_{m,\overline{\jmath}+1}}\qquad \mbox{inside} \qquad (\mathbb{C} \setminus \Delta_m)^{\prime }.
\end{equation}
The function of the right hand of (\ref{CC6}) is meromorphic in $\mathbb{C} \setminus \Delta_m$. Its zeros correspond with those of  $t_{\overline{\jmath}}$ (multiplicity included) and its poles are the zeros of $t_m$ with order equal to the multiplicity of the zero. Using the argument principle, from (\ref{CC6}) it follows that for each $j\in \{0,\ldots,m-1\}$ if $\zeta$ is
a zero of $t_j$ of multiplicity $\kappa$ then $\zeta$ attracts exactly $\kappa$ zeros of $a_{{\bf n},m}$ as $|{\bf n}| \to \infty, {\bf n} \in \Lambda(\overline{\jmath}),$ and the remaining zeros of $a_{{\bf n},m}$ accumulate of $\Delta_m$. This is true for each $\overline{\jmath}$. Hence the statement about the zeros of $a_{{\bf n},m}$ is valid for ${\bf n} \in \Lambda$ and (\ref{cu}) is satisfied.

Combining (\ref{cu}), the knowledge we have about the asymptotic behavior of the zeros of $a_{{\bf n},m}$ and the argument principle we obtain the statement about the  asymptotic behavior of the  zeros of the $a_{{\bf n},j}, j=0,\ldots,m-1$. \hfill $\Box$
\medskip

{\bf Proof of Theorem \ref{teoAT}.} There is no loss of generality if we consider multi-indices of the form
$(n_0 -d_0,\ldots,n_m-d_m)\in \mathbb{Z}_+^{m+1}$ where $d_j = \deg t_j, j=0,\ldots,m$. We will reason by contradiction.

Let us assume that their exists an infinite sequence of multi-indices $\Lambda^{\prime} \subset \Lambda$
such that for each ${\bf n} \in \Lambda^{\prime}$
their exist polynomials $p_{{\bf n},0},\ldots,p_{{\bf n},m},\deg p_{{\bf n},j} \leq n_j -d_j -1$
with real coefficients, not all identically equal to zero,  for which $ p_{{{\bf n}},0} t_0 +\sum_{j=1}^m p_{{{\bf n}},j}   t_j \hat{s}_{1,j}$  has at least $|{\bf n}|  - D$ zeros on a certain interval $\Delta \subset \mathbb{R}\setminus \Delta_{1}$, where $D = \sum_{j=0}^m d_j$. Let $w_{\bf n}$ be the polynomial whose zeros are those of $p_{{{\bf n}},0} t_0 +\sum_{j=1}^m p_{{{\bf n}},j}   t_j \hat{s}_{1,j}$ on $\Delta$. Then, it is easy to check that
\[ (p_{{{\bf n}},0} t_0 +\sum_{j=1}^m p_{{{\bf n}},j}   t_j \hat{s}_{1,j})/w_{\bf n} = {\mathcal{O}}(1/z^{|{{\bf n}}|- N_{\bf n}-D +1}) \in \mathcal{H}(\mathbb{C} \setminus \Delta_1).
\]
Therefore, the polynomials $p_{{\bf n},0}, \ldots, p_{{\bf n},m}$ fulfill a')-b'), with an extra power of $1/z$ in the right hand of b').

Let $\overline{\jmath}$ and $\Lambda(\overline{\jmath})$ be defined as in the proof of Theorem \ref{CTI}. Obviously,
$\Lambda(\overline{\jmath})$ must contain infinitely many multi-indices in $\Lambda^{\prime}$ for some
$\overline{\jmath} \in \{0,\ldots,m\}$. Fix
$\overline{\jmath}$ so that this occurs. Arguing as in the proof of Theorem  \ref{CTI}, we obtain that $p_{{\bf n},\overline{\jmath}}t_{\overline{\jmath}}$  has at least $n_{\overline{\jmath}} - D$ zeros on $\Delta_m$ and for all sufficiently large $|{\bf n}|$ as many zeros close to each one of the zeros  of $t_j, j=0,\ldots,m$, as their multiplicity. Therefore, for all sufficiently large $|{\bf n}|, {\bf n} \in \Lambda(\overline{\jmath}) \cap \Lambda^{\prime}$ we have that $\deg p_{{\bf n},\overline{\jmath}}t_{\overline{\jmath}} = n_{\overline{\jmath}}$. This contradicts the fact that by construction $\deg p_{{\bf n},\overline{\jmath}}t_{\overline{\jmath}} \leq n_{\overline{\jmath}}-1$. Thus our initial assumption is false and the statement of the theorem  true. \hfill $\Box$

\begin{remark} Theorem \ref{TCTII} can be extended to multipoint type II Hermite-Pad\'e approximation of ${\bf f}= \widehat{\bf s} + {\bf r}$.
\end{remark}

\end{document}